\documentclass[reqno,11pt,a4]{amsart}
\usepackage{amsthm,amsfonts,amssymb,amsmath,latexsym,marginnote,mathtools}
\mathtoolsset{showonlyrefs=true}
\usepackage[dvips]{graphicx}
\usepackage{hyperref}
\usepackage[margin=2.5cm]{geometry}
\usepackage{enumerate}
\usepackage{color}
\usepackage{comment}

\newtheorem{theorem}{Theorem}[section]

\newtheorem{proposition}[theorem]{Proposition}

\newtheorem{remark}[theorem]{Remark}

\def\R {\mathbb{R}}

\def\N {\mathbb{N}}

\title[Multiple-wise interactions for multi-agent systems]{Multi-agent systems with multiple-wise interaction: Propagation of chaos and macroscopic limit}

\author{Thierry Paul}
\address{CNRS Laboratoire Ypatia des Sciences Mathematiques LYSM, Rome}
\email{thierry.paul@sorbonne-universite.fr}

\author{Stefano Rossi}
\address{ETH Zurich}
\email{stefano.rossi@math.ethz.ch}

\author{Emmanuel Trélat}
\address{Sorbonne Université, Paris}
\email{emmanuel.trelat@sorbonne-universite.fr}

\begin{document}
\begin{abstract}
 We consider interacting multi-agent systems where the interaction is not only pairwise but involves simultaneous interactions among multiple agents (multiple-wise interaction).
 By passing through the mesoscopic and macroscopic limits with a fixed multiple-wise interaction of order $m$, we derive a macroscopic equation in the limit $m\to\infty$, capturing the dominant effects in large-size multiple-wise order. 
\end{abstract}
\date{\today}

\maketitle
\tableofcontents

\section{Introduction}

Many mathematical models for systems of particles or agents, as in physics or biology, typically assume pairwise interactions among the agents. In this article, we investigate more general models involving what we entitle \emph{multiple-wise} interactions.

Specifically, given $N\in \mathbb{N}^*$ and $m \in \{1,\dots,N\}$, we consider a system of agents $\{\xi^N_i(t)\}_{i=1}^N \subset \R^d$ (for some $d \in \mathbb{N}^*$) communicating via multi-body interactions (of order $m$) of the type
\begin{equation}
\label{sec0:maineq}
	\dot\xi_i^N(t)=\frac{1}{N^m}\sum_{j_1, \dots, j_m=1}^N G^N_{i,j_1,\dots,j_m}\left(t,\xi_i^N(t),\xi^N_{j_1}(t), \dots,\xi^N_{ j_m}(t)\right), \quad i \in \{1,\dots,N\},
\end{equation}
where
$
G^N_{i,j_1,\dots,j_m}: \mathbb{R}_t\times \mathbb{R}_\xi^{d(m+1)} \to \mathbb{R}^d
$
represents the multiple-wise interaction between the agent $i$ and the agents $j_1,\dots,j_m$, and where the indices possibly coincide.\\

In applied sciences, such as opinion dynamics, biophysics, and generally the physics of complex systems, models incorporating multiple-wise interactions like the one above have recently gained attention (see, e.g., \cite{BATTISTON20201}). These models provide new deep insights into the emergent phenomena related to complex systems (see, e.g., \cite{RUSSO24}).

Even the case $m=1$,  corresponding to pairwise interactions, is of significant interest in fields like opinion dynamics. 
Here, the value $\xi^N_i(t)$ represents the opinion of agent $i$, which naturally depends on its label. 
However, group behaviors often cannot be fully explained by pairwise interactions alone, since opinions are frequently governed by multiple-wise interactions
($m>1$) of the type described above.\\

This paper focuses on the statistical properties of systems as described by \eqref{sec0:maineq}, particularly in the large $N$ and $m$ limits. For large $N$, the statistical behavior of the system does not depend on the exact value of $N$, enabling reduced descriptions via partial differential equations as $N \to +\infty$.

In the pairwise case ($m=1$) these limits have been extensively studied in the literature, both for the classical undistinguishable case (see, e.g., \cite{Dobr}), and in recent works for the case of distinguishable particles (see, e.g., \cite{PT24}). In this work, we extend the analysis to the multiple-wise case ($m > 1$) with distinguishable particles, exploring whether a limit equation can be derived as $m \to +\infty$, preserving the key characteristics of the interaction type.

To achieve this, we consider three distinct limits: the mesoscopic limit, the macroscopic limit, and a large multiple-wise order limit where, after taking $N \to +\infty$, we also let $m \to +\infty$. 

At the mesoscopic level, agent labels and opinions are treated as independent variables distributed according to a statistical distribution on the label-opinion phase space. By taking marginals of the symmetrization (see eq. \eqref{sec1:symm}) of the distribution function, we show that, as $N \to +\infty$, the discrete dynamics is well-approximated by an $m$-dependent multiple-wise interacting Vlasov-type equation.

At the macroscopic level, opinions are functions of agent labels, distributed according to a probability density. We show how to pass from the mesoscopic to a macroscopic description, yielding an $m$-dependent integro-differential equation for the opinion functions. 

Finally, as $m \to +\infty$, we derive a limiting differential equation for the limit opinion function influenced by infinitely many other agents. This equation will retain the mean features of the interactions and label distributions in an asymptotic vector field.

By combining the three asymptotic procedures, we obtain the main result of the work (see part \eqref{sec0:thm2pt4} of Theorem \ref{sec0:thm2}), which is a convergence theorem from the single-particle distribution function related to the dynamics of $N$ particles with $m$-wise interaction, to the limit equation satisfied by the opinion function in the limit as $N$ and $m \to +\infty$.\\

We introduce the models and some notation in the following Subsections \ref{sec0:sub1} and \ref{sec0:sub2}. In Section \ref{sec2:state} we present the main results given by Theorems \ref{sec0:thm1} and \ref{sec0:thm2}. The remaining Sections \ref{sec3:sec} and \ref{sec4:sec} of the paper are dedicated to the proofs.

\subsection{Models and assumptions}
\label{sec0:sub1}
In the following, for any $m \in \{1,\dots,N\}$ we assume that there exist a complete metric space $(\Omega, \text{d}_\Omega)$ and a continuous function
\begin{equation}
    \label{sec0:mainass0}
G^{(m)}: \R_t \times (\Omega_x \times \R_\xi^{d})^{m+1} \to \R^d
\end{equation}
that is globally Lipschitz on $\Omega_x^{m+1}\times\R_\xi^{d(m+1)}$ uniformly in every compact subset of $\R_t$ and such that for every $N$, there exist $x_1,\dots x_N \in \Omega$ such that
\begin{equation}
    \label{sec0:mainass}
G^{(m)}(t,x_{i}, \xi,x_{j_1},\xi_1, \dots, x_{j_m},\xi_{m})=G^N_{i,j_1,\dots,j_m}(t,\xi,\xi_1\dots, \xi_{m}),
\end{equation}
for any $t \in \R$, $(\xi,\xi_1,\dots,\xi_{m})\in \R^{d(m+1)}$ and any $\{i,j_1,\dots,j_m\}\subset \{1,\dots,N\}$ and where $G^N_{i,j_1,\dots,j_m}$ is defined as in \eqref{sec0:maineq}. 

We suppose that:
\begin{enumerate}[(A)]
\item
\label{sec0:ass1}
There exists $M>0$ independent of $m$ such that
\begin{equation}
\label{sec0:assM}
\left\|G^{(m)}\right\|_{L^\infty(\R_t\times(\Omega_x \times \R_\xi^d)^{m+1})}\le M;
\end{equation}
\item 
\label{sec0:ass2}
$G^{(m)}$ is globally Lipschitz continuous on $\left(\Omega_x \times \R^d_\xi\right)^{m+1}$ with a Lipschitz constant $\text{Lip}_{x,\xi}(G)$ independent of $m$.
\end{enumerate}

In particular $G^{(m)}$
is globally Lipschitz on $\R_\xi^{d(m+1)}$, with respect to the norm
\begin{equation}
    \label{sec0:lq}
\|\Xi_{m+1}\|_{\ell_q}:=\|(|\xi_1|, \dots, |\xi_{m+1}|)  \|_{\ell_q},\quad \Xi_{m+1}:=(\xi_1,\dots,\xi_{m+1})\in \mathbb{R}^{d(m+1)}, \quad  q \in [1,\infty].
\end{equation}

We pass to describe the three models at the micro/meso/macroscopic level.\\

\textbf{Microscopic level:}
By \eqref{sec0:mainass}, we can equivalently define the multi-agent system \eqref{sec0:maineq} as
\begin{equation}
\begin{aligned}
\label{sec0:withlabel}
    \dot x^N_i(t)&=0 \\
    \dot \xi^N_i(t)&=\frac{1}{N^m}\sum_{j_1, \dots, j_m=1}^N G^{(m)}\left(t,x^N_i,\xi^N_i(t),x^N_{j_1},\xi^N_{j_1}(t),\dots,x^N_{j_m},\xi^N_{ j_m}(t)\right), \quad i \in \{1,\dots,N\}.
\end{aligned}
\end{equation}

In this way the (stationary) labels $\{x^N_i\}_{i=1}^N$ related to the opinions $\{\xi_i^N\}_{i=1}^N$ illustrate the distinguishability of the particles.

Notice that, since the interaction $G^{(m)}$
  is assumed to be globally Lipschitz, 
  the classical Cauchy-Lipschitz theory ensures that the solutions to \eqref{sec0:withlabel} are globally defined in time.\\
  
\textbf{Mesoscopic level:}
Given $t \in \R$, we introduce the empirical measure at time $t$ associated with system \eqref{sec0:withlabel}; this is given by
\begin{equation}
\label{sec0:empirical}
\mu_N(t,x,\xi):=\frac{1}{N} \sum_{i=1}^N \delta_{x^N_i}(x)\otimes\delta_{\xi^N_i(t)}(\xi) \in \mathcal{P}(\Omega_x \times \R^d_\xi).
\end{equation}
For a regular function $\varphi: \Omega_x\times\mathbb{R}^d_\xi \to \mathbb{R}$, we have that
\begin{equation}
    \begin{aligned}
\frac{d}{dt} \left \langle \varphi, \mu_N(t) \right \rangle&=\frac{1}{N^{m+1}}\sum_{i,j_1, \dots, j_m=1}^N G^{(m)}\left(t,x^N_i,\xi^N_i(t),x^N_{j_1},\xi^N_{j_1}(t), \dots,x^N_{j_m},\xi^N_{ j_m}(t)\right) \cdot \nabla_{\xi} \varphi(x_i^N,\xi^N_i(t))\\
&=\left \langle \mathcal{X}^{[m]}[\mu_N(t)] \cdot\nabla_{\xi}\varphi, \mu_N(t)\right \rangle,
\end{aligned}
\end{equation}
where $\mathcal{X}^{[m]}[\mu_N]$ is the multiple-wise mean-field force of order $m$ given by
\begin{equation}
\label{vf}
\mathcal{X}^{[m]}[\mu_N](t,x,\xi):= \int_{(\Omega_x \times\mathbb{R}_\xi^{d})^m} G^{(m)}(t,x,\xi,x_1,\xi_1,\dots,x_m,\xi_m) d \mu^{\otimes m}_N(x_1,\xi_1,\dots,x_m,\xi_m),
\end{equation}
so that $\mu_N(t)$ weakly solves the following mesoscopic equation of Vlasov-type:
\begin{equation}
\label{sec0:meanfield}
\partial_t f(t,x,\xi) + \nabla_{\xi} \cdot \left(\mathcal{X}^{[m]}[f(t)](t,x,\xi) f(t,x,\xi)\right)=0.
\end{equation}

In the following we want to justify the use of equation \eqref{sec0:meanfield} by proving the existence and uniqueness of its weak measure-valued solutions and showing the following mesoscopic limit: If $\mu_N(0)$ is converging to a given probability distribution $f_0 \in \mathcal{P}(\Omega_x \times \R_\xi^d)$ as $N \to +\infty$, the same holds for a time $t>0$, i.e. $\mu_N(t)$ converges to the unique solution $f(t) \in \mathcal{P}(\Omega_x \times \R_\xi^d)$ to \eqref{sec0:meanfield} with initial datum $f_0$. 
This can be proved quantitatively by obtaining a Dobrushin-type estimate (as referred by the seminal paper \cite{Dobr}).\\

A different (but related) point view is the one of propagation of chaos:
Given $f(0)\in \mathcal{P}(\Omega_x\times \R^d_\xi)$, 
assume that at the initial time the agents are
independently and identically distributed with law $F^N(0):=f(0)^{\otimes N} \in \mathcal{P}((\Omega_x \times \R_\xi^d)^N)$.
For a positive time $t>0$, the statistical distribution of the agents is given by $F^N(t)\in\mathcal{P}((\Omega_x \times \R_\xi^d)^N)$ weak solution of the following Liouville equation
\begin{equation}
\label{sec0:Liouville}
    \partial_t F^N(t,Z_N) + \sum_{i=1}^N \frac{1}{N^m} \sum_{j_1,\dots,j_m=1}^N \nabla_{\xi_i}\cdot \left(G^{(m)}(t,x_i,\xi_i,x_{j_i},\xi_{j_1},\dots,x_{j_m},\xi_{j_m}) F^N(t,Z_N)\right)=0,
\end{equation}
where, for $k \in \{1,\dots,N\}$, $Z_k:=(x_{1},\xi_1,\dots,x_k,\xi_k)$.
In general, for positive times $t>0$, the distribution function is no longer symmetric since the particles are distinguishable and the function $G^{(m)}$ is not symmetric in the couples $(x_i,\xi_i)_{i=1}^N$. For this reason we introduce the symmetrization of the distribution function given by
\begin{equation}
    \label{sec1:symm}
    F^N(t)^s(Z_N):= \frac{1}{N!} \sum_{\sigma \in \mathfrak{S}_N} F^N(t,\sigma(Z_N))
\end{equation}
where $\mathfrak{S}_N$ is the set of permutations on $\{1,\dots,N\}$ and $\sigma(Z_N)=(x_{\sigma(1)},\xi_{\sigma(1)},\dots,x_{\sigma(N)}, \xi_{\sigma(N)})$.
Defining
\[
F^N(t)^s_{N:k}(Z_k):=\int_{(\Omega_x\times \R_\xi^d)^{N-k}} F^N(t)^s(Z_N) dx_{k+1}d\xi_{k+1}\dots dx_Nd\xi_{N},
\]
we expect that, for $N$ large, the behavior of the small group of $k$ particles will become nearly independent, so that at time $t$ and as $N \to +\infty$, each particle will follow the same statistical distribution $f(t)$ given by the solution of \eqref{sec0:meanfield} with initial datum $f(0)$: i.e., in a suitable topology, $F^N(t)^s_{N:k}$ will converge to $f^{\otimes k}(t)$. 
All of this will be rigorously proved in Theorem \ref{sec0:thm1}.\\

\noindent\textbf{Macroscopic level:}
Let us consider again the mesoscopic equation \eqref{sec0:meanfield}.
By taking a monokinetic weak solution to \eqref{sec0:meanfield} of the form
\begin{equation}
    f_{\text{mon}}(t,x,\xi):=\nu(x) \otimes \delta_{y(t,x)}(\xi)\in \mathcal{P}(\Omega_x\times \R_\xi^d),
\end{equation}
with $\nu \in \mathcal{P}(\Omega_x)$ and $y(t,x)$ a measurable function of class $\mathcal{C}^1$ respect to time, it is not difficult to show that $y(t,x)$ solves the macroscopic equation
\begin{equation}
\label{sec2:opinion}
\partial_t y(t,x) = \int_{\Omega^m} G^{(m)}(t,x,y(t,x),x_1,y(t,x_1),\dots,x_m,y(t,x_m)) d\nu^{\otimes m}(x_1,\dots,x_m).
\end{equation}
In Theorem \ref{sec0:thm2}, after proving the existence and uniqueness of solutions to \eqref{sec2:opinion}, we study the limit as the order of the multiple-wise interaction 
$m$ tends to infinity, obtaining a limit equation which retains the features of the interaction and of the distribution of the labels (see \eqref{sec1:limitm} in part \eqref{sec0:thm2pt4} in Theorem \ref{sec0:thm2}).

\subsection{Notation}
\label{sec0:sub2}

Given a measurable space $X$, we denote by $\mathcal{P}(X)$ the set of probability measures on $X$. In order to quantify the difference between two probability measures $\mu, \nu \in \mathcal{P}(\Omega_x \times \R_\xi^d)$, we introduce the Wasserstein distance
\begin{equation}
    \label{sec0:wassdist}
    \mathcal{W}_p\left(\mu,\nu\right):=\inf_{\pi \in \mathcal{C(\mu,\nu)}}\left(\int_{(\Omega_x \times \R_\xi^d)^2} \left(\text{d}_\Omega(x^1-x^2) + |\xi^1-\xi^2|\right)^pd\pi(x^1,\xi^1,x^2,\xi^2)\right)^{\frac1p},
\end{equation}
where $\mathcal{C}(\mu,\nu)$ is the set of all couplings between $\mu$ and $\nu$: i.e. probability measures on the product space $(\Omega_x \times \R_\xi^d)^2$ with first and second marginal respectively $\mu$ and $\nu$.

We will also compare two probability measures $\mu_k,\nu_k\in \mathcal{P}\left( (\Omega_x\times \R^d_\xi)^k\right)$: In this case,
denoting $X_k^{i}:=(x_1^i,\dots,x_k^{i})\in\Omega_x^k$, $\Xi_k^i:=(\xi_1^i,\dots,\xi_k^i)\in \R_\xi^{dk}$ and $Z_k^{i}:=(x_1^i,\xi_1^i,\dots,x_k^{i},\xi_k^i)\in(\Omega_x\times \R^d_\xi)^k$ for $i=1,2$,
we define the $q$-modified Wasserstein distance
\begin{equation}
    \label{sec0:wassdistmod}
    \mathcal{W}^{[q]}_p\left(\mu_k,\nu_k\right):=\inf_{\pi \in \mathcal{C}(\mu_k,\nu_k)}\left(\int_{(\Omega_x \times \R_\xi^d)^{2k}} \left(\text{d}^{(q)}_\Omega(X_k^1-X_k^2) + \|\Xi_k^1-\Xi_k^2\|_{\ell^q}\right)^pd\pi(Z_k^1,Z_k^2)\right)^{\frac1p},
\end{equation}
where, given $q\ge1$, the $\ell^q$ norm is defined in \eqref{sec0:lq}, while
\[
\text{d}^{(q)}_\Omega(X_k^1-X_k^2):=\left \|\left(\text{d}_\Omega(x^1_1-x^2_1),\dots,\text{d}_\Omega(x_k^1-x^2_k) \right) \right\|_{\ell^q}.
\]

Given $f\in \mathcal{C}^0(\mathbb{R}; \mathcal{P}(\Omega_x\times\R_\xi^d))$ and a globally Lipschitz vector field $\mathcal{X}[f(t)]$ as in \eqref{vf}, we also introduce the flow $\varphi_{f_0}(t): \Omega_x \times \R^d_\xi\to\Omega_x\times \R^d_\xi$ for $t \in \R$, given by $\varphi_{f_0}(t)(x_0,\xi_0):=(x_0,\xi(t))$ where 
\begin{equation}
\label{sec0:flowmap}
    \begin{cases}
        \dot x(t)=0\\
        \dot \xi(t)=\mathcal{X}^{[m]}[f(t)](t,x(t),\xi(t))\\
        x(0)=x_0,\, \xi(0)=\xi_0, \quad (x_0,\xi_0)\in\Omega_x \times \R^d_\xi.
    \end{cases}
\end{equation}
Moreover, if $(X,\mu)$ is a probability space, $Y$ a measurable space and $\phi:X \to Y$ a measurable function, we define the push-forward $\phi_*\mu \in \mathcal{P}(Y)$ of $\mu$ along $\phi$ by
\[
\phi_*\mu(A):=\mu(\phi^{-1}(A)), \quad A \subset Y\quad \text{measurable}.
\]

\section{Statements of the main results}
\label{sec2:state}

We now state the main results of our paper. First, we validate the mesoscopic equation by deriving it from the multiple-wise multi-agent system in \eqref{sec0:maineq} and by proving the existence and uniqueness of its weak solutions.

\begin{theorem}[Existence and uniqueness, Dobrushin estimate and propagation of chaos]
\label{sec0:thm1}
Given $m \in \mathbb{N}^*$, let $G^{(m)}$ be as in \eqref{sec0:mainass0} satisfying assumptions \eqref{sec0:ass1} and \eqref{sec0:ass2}.
The following results hold:

\begin{enumerate}[(i)]
\item 
\label{sec0:thm1pt1}
(Existence and uniqueness).
Given $f_0(x,\xi) \in \mathcal{P}(\Omega_x\times\mathbb{R}_\xi^d)$, there exists a unique global solution $f \in \mathcal{C}([0,+\infty); \mathcal{P}(\Omega_x \times \mathbb{R}_\xi^d))$ to the Vlasov-type equation \eqref{sec0:meanfield} with initial datum $f_0$. Equivalently, $f$ is the unique solution of the equation:
\begin{equation}
\label{sec1:mono}
f(t)=\varphi_{f_0}(t)* f_0, \quad t\in\R
\end{equation}
where $\varphi_{f_0}(t)$ is the Vlasov flow defined as in \eqref{sec0:flowmap}.
\item
\label{sec0:thm1pt2}
(Dobrushin-type estimate).
    Given $f(0) \in \mathcal{P}(\Omega_x \times \R^d_x)$, let $f(t,x,\xi)\in\mathcal{P}(\Omega_x \times \R^d_\xi)$ be the weak solution of the Vlasov-type equation \eqref{sec0:meanfield} with initial datum $f(0)$ and $\mu_N(t,x,\xi)$ be the empirical measure defined in \eqref{sec0:empirical} with initial datum $\mu_N(0)$.
Then, for $m \in \{1, \dots, N\}$ and $p \in [1,+\infty), q\in[1,+\infty]$ , it holds that
\begin{equation}
\label{mfl}
\mathcal{W}_p(\mu_N(t), f(t)) \le C e^{C \text{Lip}_{x,\xi}(G) (1+m^{1/q})t} \mathcal{W}_p(\mu_N(0), f(0)).
\end{equation}

\item
\label{sec0:thm1pt3}
(Propagation of chaos).
Assume that $\Omega_x \subset \R^d$ and, for $k \in \{1, \dots, N\}$, let $F^N(t)^s_{N:k}$ be the $k$-marginal of the symmetrization of the solution to the Liouville equation \eqref{sec0:Liouville} with initial datum $F^N_0(Z_N)=\prod_{i=1}^Nf_0(x_i,\xi_i)$, with $f_0 \in \mathcal{P}(\Omega_x \times \R^d_\xi)$.
Let $f(t)\in\mathcal{P}(\Omega_x \times \R^d_\xi)$ be the unique solution to \eqref{sec0:meanfield} with initial datum $f_0$, then, for $p \in [1,+\infty)$, $q \in [1,+\infty]$ and $z>0$,
\begin{equation}
\label{sec0:result2}
\begin{aligned}
\mathcal{W}^{[q]}_p&(F^N(t)^s_{N:k}, f(t)^{\otimes k})\le
C(km)^{\frac{1}{q}} \text{Lip}_{x,\xi}(G) e^{\text{Lip}_{x,\xi}(G)(1+m^{1/q}) t}
\left(\int_0^t \mathcal{M}_z(f(\tau))^{1/z} d\tau\right)\times\\
&\times
\begin{cases}
N^{-\frac{1}{2p}}+ N^{-(z-p)/(pz)} &\, \text{if}\quad p>\frac{d}{2} \quad \text{and} \quad z\neq 2p\\
N^{-\frac{1}{2p}}(\log(1+N))^{1/p}+ N^{-(z-p)/(pz)} &\,\text{if} \quad p=\frac{d}{2} \quad \text{and} \quad z\neq 2p \\
N^{-\frac{1}{d}}+ N^{-(z-p)/(pz)} &\,\text{if}\quad p \in (0, \frac{d}{2})\,\, \text{and}\,\, z\neq d/(d-p),
\end{cases}
    \end{aligned}
\end{equation}
where $\mathcal{M}_z(f):= \int |x|^z df(x)$.
\end{enumerate}
\end{theorem}

\begin{remark}
    Theorem \ref{sec0:thm1} can be easily extended to a locally Lipschitz vector field $G$, by analyzing the time of existence of the $N$-particle system and controlling the temporal evolution of the supports of the measures (see \cite{PT24} for the case \(m=1\)). This extension is especially relevant when studying second-order systems with interactions depending on the velocities of the agents.
\end{remark}
\begin{remark}
       Notice that in \eqref{sec0:result2}, we can also consider the joint limit in $m$ and $N$ when the interaction order \(m\) doesn't grow too fast in \(N\), specifically with
       \begin{equation}
       \label{sec2:m(N)}
       m \le m(N) = \lfloor (\text{Lip}_{x,\xi}(G)T)^{-1} \log(N^\alpha) \rfloor^{q}
       \end{equation}
       for \(\alpha < \frac{1}{2}\) (see part \eqref{sec0:thm2pt4} of Theorem \ref{sec0:thm2} for further details).
\end{remark}

We now state the main result concerning the macroscopic equation \eqref{sec2:opinion}. In particular, parts \eqref{sec0:thm2pt3} and \eqref{sec0:thm2pt4} of Theorem \ref{sec0:thm2} deal with the large multiple-wise order limit.

\begin{theorem}
\label{sec0:thm2}
    Given $m \in \mathbb{N}^*$, let $G^{(m)}$ be as in \eqref{sec0:mainass0} satisfying assumptions \eqref{sec0:ass1} and \eqref{sec0:ass2}.
        The following results hold:
    \begin{enumerate}[(i)]
        \item \label{sec0:thm2pt1}
Given $T>0$, let $\nu \in \mathcal{P}(\Omega_x)$ and $y(t,\cdot) \in \mathcal{C}^1\left([0,T];L^\infty(\Omega_x;\R^d)\right)$. Then 
$$
f_{\text{mon}}(t,x,\xi):=\nu(x)\otimes \delta_{y(t,x)}(\xi)\in\mathcal{P}(\Omega_x\times\R^d_\xi)
$$ 
is a (weak) solution of the Vlasov-type equation \eqref{sec0:meanfield} if and only if $y(t,x)$ is solution of the macroscopic equation \eqref{sec2:opinion}, $\nu_x$--almost everywhere on $[0,T]$.
   
    \item \label{sec0:thm2pt2}
Given $\nu \in \mathcal{P}(\Omega_x)$ and $y_0 \in L^\infty(\Omega_x;\mathbb{R}^d)$, there exists a unique global solution $ y^{(m)}(t,\cdot) \in L^\infty(\Omega_x;\mathbb{R}^d)$ which is continuously differentiable in time of the macroscopic equation \eqref{sec2:opinion} such that $y(0)=y_0$.

 \item 
 \label{sec0:thm2pt3}
 Given $\nu \in \mathcal{P}(\Omega_x)$ and a Lipschitz datum $y_0 \in L^\infty(\Omega_x;\R^d)$, there exists $y^\infty\in \mathcal{C}^1([0,T];L^\infty(\Omega_x;\R^d))$ and $G^\infty_\nu:[0,T]\times \Omega_x\times L^\infty(\Omega_x)\to \R^d$
    such that
    \begin{equation}
    \label{sec1:limitm}
        \partial_ty^\infty(t,x)=G^\infty_\nu[y^\infty(t)](t,x)
    \end{equation}
    and, up to a subsequence,
    \begin{equation}
        \label{sec3:multilim}
        \lim_{m \to \infty}\sup_{t \in [0,T]}\|y^{(m)}(t,\cdot)-y^\infty(t,\cdot)\|_{L^\infty(\Omega)}=0,
    \end{equation}
    where $y^{(m)} \in \mathcal{C}^1([0,T],L^\infty(\Omega_x;\R^d))$ is the unique solution to the macroscopic equation \eqref{sec2:opinion} with multiple-wise order $m$.
    \item 
    \label{sec0:thm2pt4}
    Let $\Omega_x \subset \R^d$, given $\nu \in \mathcal{P}(\Omega_x)$ and $y_0 \in L^\infty(\Omega_x)$, let $F^{N,(m)}(t)^s_{N:1}$ be the first marginal of the symmetrization of the solution to the Liouville equation \eqref{sec0:Liouville} with multiple-wise order $m$ and monokinetic initial datum
    \[
    F^N_{\text{mon}}(0,X_N,\Xi_N):=\prod_{i=1}^N \nu(x_i)\otimes \delta_{y_0(x_i)}(\xi_i).
    \]
    Moreover, let us denote by $Y_{N:1}^{(m)}(t,x):=\int\xi F^{N,(m)}(t)^s_{N:1}(t,x,\xi) d\xi$ its first moment in $\xi$. Then, for any $\phi \in \mathcal{C}(\Omega_x)$
    and $\alpha < \frac12$,
    $$
    \lim_{\substack{N,m \to + \infty\\  m\leq m(N)}}
    \sup_{t \in [0,T]}\int_{\Omega_x} \phi(x)\left(Y_{N:1}^{m}(t,dx)-y^\infty(t,x)\nu(dx) \right)=0,
    $$
    where $m(N)$ is defined in \eqref{sec2:m(N)}.
    \end{enumerate}
   
\end{theorem}
\begin{remark}
We emphasize that -- although the solutions of the mesoscopic Vlasov equation \eqref{sec0:meanfield} with measure data exist globally in time, thanks to part \eqref{sec0:thm1pt1} of Theorem \ref{sec0:thm1} -- the preservation of the monokinetic nature of solutions is not guaranteed globally but only up to a time \( T > 0 \). In particular, the statement of part \eqref{sec0:thm2pt1} of Theorem \ref{sec0:thm2} holds as long as there is no blow-up of the monokinetic property.
\end{remark}

\section{Proof of Theorem \ref{sec0:thm1}: Validation of the mesoscopic equation}
\label{sec3:sec}

\textbf{Proof of part \eqref{sec0:thm1pt1} of Theorem \ref{sec0:thm1}: existence and uniqueness of weak solutions.}
The existence and uniqueness of the solution follow by a classical contraction argument by defining the map
\begin{equation}
    \mathcal{F}: L^\infty([0,+\infty); \mathcal{P}(\Omega_x\times \mathbb{R}_\xi^d)) 
    \to
    L^\infty([0,+\infty); \mathcal{P}(\Omega_x\times \mathbb{R}_\xi^d))
\end{equation}
such that $\mathcal{F}(\nu)(t):=\varphi_{\nu(0)}(t)_*\mu_0$, for $\nu \in \mathcal{C}([0,+\infty); \mathcal{P}(\Omega_x\times\mathbb{R}_\xi^d))$. The contractivity and the global existence follow by the assumption of global Lipschitzianity of $G^{(m)}$.\qed\\

\textbf{Proof of part \eqref{sec0:thm1pt2} of Theorem \ref{sec0:thm1}: Dobrushin-type estimate.}

Let $\Pi_0 \in \mathcal{C}(\mu_N(0),f(0))\in \mathcal{P}((\Omega_x\times \R^d_\xi)^2)$ be a general coupling between $\mu_N(0)$ and $f(0)$
and let $\Phi^N(t)$ be the particle flow at time $t$ generated by the particle system in \eqref{sec0:maineq} and $\varphi_{f_0}(t)$ the Vlasov flow given in \eqref{sec1:mono}.
Then
the probability measure on $\mathcal{P}((\Omega_x\times \R_\xi^d)^2)$ given by 
\[
\Pi_t(x^1,\xi^1,x^2,\xi^2):=\left( \Phi_t^N \otimes \varphi_{f_0}(t)\right)_* \Pi_0(x^1,\xi^1,x^2,\xi^2),
\]
is a coupling between $\mu_N(t)$ and $f(t)$, solutions to the mesoscopic Vlasov-type equation \eqref{sec0:meanfield} with initial data rescpectively $\mu_N(0)$ and $f(0)$.

Indeed notice that the marginals of $\Pi_t$ on the first $(x^1,\xi^1)$ and second $(x^2,\xi^2)$ couple of variables are both solutions of the mesoscopic Vlasov equation \eqref{sec0:meanfield}. By the uniqueness in part \eqref{sec0:thm1pt1} of Theorem \ref{sec0:thm1}, it follows that $\Pi_t$ is a coupling between $\mu_N(t)$ and $f(t)$, i.e. $\Pi_t\in \mathcal{C}(\mu_N(t),f(t))$ (see also \cite[Lemma 3.2]{GMP16}).
 
By definition of Wasserstein distance, we have that
$$
\mathcal{W}_p(\mu_N(t), f(t))^p \le \int_{(\Omega \times \mathbb{R}^d)^2} \left(\text{d}_\Omega(x^1-x^2)+|\xi^1-\xi^2|\right)^p d \Pi_t(x^1,\xi^1,x^2,\xi^2)
$$
and, by definition of push-forward,
$$
    \mathcal{W}_p(\mu_N(t), f(t))^p \le \int_{(\Omega \times \mathbb{R}^d)^2}\left(\text{d}_\Omega(x^1-x^2)+|\xi^1(t)-\xi^2(t)| \right)^p d\Pi_0(x^1,\xi^1,x^2, \xi^2),
$$
where $\xi^1(t):=\pi_\xi(\Phi^N_t(x^1,\xi^1))$, $\xi^2(t):=\pi_\xi(\varphi_{f(0)}(t)(x^2,\xi^2))$ and $\pi_\xi:\Omega_x\times \R^d_\xi\to \R^d_\xi$ is the projection on the second component. We have
\begin{align}
\left|\xi^1(t) - \xi^2(t)\right|\le \left|\xi^1 - \xi^2\right| 
&+ \int_0^t \Bigl(\left|\mathcal{X}^{[m]}[\mu_N(\tau)](\tau, x^1,\xi^1(\tau)) - \mathcal{X}^{[m]}[f(\tau)](\tau,x^1, \xi^1(\tau))\right| \\
&+ \left|\mathcal{X}^{[m]}[f(\tau)](\tau, x^1,\xi^1(\tau)) - \mathcal{X}^{[m]}[f(\tau)](\tau, x^2,\xi^2(\tau))\right| \Bigr)d \tau.
\end{align}
For $q \in [1,+\infty]$, by using the regularity of $G^{(m)}$, the first term in the integral can be bounded by
\begin{align*}
 \int_{\left(\Omega_x\times \mathbb{R}^d_\xi\right)^m} 
&|G^{(m)}(\tau,x^1,\xi^1(\tau),Z_m) - G^{(m)}(\tau,x^1, \xi^1(\tau),Z_m')|
d \tilde\Pi(Z_m,Z'_m)\\
&\le \text{Lip}_{x,\xi}(G) \int_{(\Omega_y \times \R_\zeta^d)^m} 
\left(\text{d}^{(q)}_\Omega(X_m-X'_m)+\|\Xi_m-\Xi_m'\|_{\ell^q}\right)
d \tilde\Pi(Z_m,Z'_m),
\end{align*}
for any coupling $\tilde\Pi\in\mathcal{P}((\Omega_x \times \R^d_\xi)^{2m})$ between $\mu^{\otimes m}_N(\tau)$ and $f^{\otimes m}(\tau)$ and
where $X_m=(x_1,\dots,x_m)\in \Omega_x^m$, $\Xi_m=(\xi_1,\dots,\xi_m)\in(\R_\xi^d)^{m}$ and $Z_m=(x_1,\xi_1,\dots,x_m,\xi_m)\in (\Omega_x\times\R^{d}_\xi)^m$.

The second term can also be bounded thanks to the regularity of $G^{(m)}$:
$$
|
\mathcal{X}^{[m]}[f(\tau)](\tau, x^1,\xi^1(\tau))
 - \mathcal{X}^{[m]}[f(\tau)](\tau,x^2, \xi^2(\tau))
|
 \le \text{Lip}_{x,\xi}(G)\left(\text{d}_\Omega(x^1-x^2)+|\xi^1(\tau)- \xi^2(\tau)|\right).
$$
By convexity of the $p$-norm and the Minkowski's integral inequality, we conclude
\begin{align*}
\mathcal{W}_p(\mu_N(t), f(t)) &\le C \Bigl(\mathcal{W}_p(\mu_N(0), f(0)) + \text{Lip}_{x,\xi}(G) \int_0^t \mathcal{W}_p(\mu_N(\tau), f(\tau)) d\tau + \\
&\qquad \text{Lip}_{x,\xi}(G)\int_0^t \mathcal{W}^{[q]}_p(\mu_N^{\otimes m}(\tau), f^{\otimes m}(\tau))d\tau \Bigr)\\
&\le 
C \left(\mathcal{W}_p(\mu_N(0), f(0)) + \text{Lip}_{x,\xi}(G)(1 + m^{\frac{1}{q}}) \int_0^t \mathcal{W}_p(\mu_N(\tau), f(\tau))d\tau \right),
\end{align*}
where we used that $\mathcal{W}^{[q]}_p\left(\mu_N^{\otimes m}(\tau), f^{\otimes m}(\tau)\right) \le m^{\frac1q}\mathcal{W}_p(\mu_N(\tau), f(\tau))$ (see \cite[Lemma 14 (Appendix A.1.4)]{PT24}). We get \eqref{mfl} by the Grönwall's lemma.\qed\\

\textbf{Proof of part \eqref{sec0:thm1pt3} of Theorem \ref{sec0:thm1}: Propagation of chaos.}
Since $F^N(t)^s, f(t)^{\otimes N}$ are both symmetric probability measures belonging to $\mathcal{P}((\Omega_x \times \R^d_\xi)^N)$, i.e. $$\sigma_*F^N(t)^s=F^N(t)^s, \quad\sigma_*f(t)^{\otimes N}=f(t)^{\otimes N}$$ for every $\sigma \in \mathfrak{S}_N$,
we have (see \cite[Lemma 21 (Appendix A.2.2)]{PT24}):
\begin{equation}
\label{proof1:margest}
\mathcal{W}^{[q]}_p\left(F^{N}(t)^s_{N:k}, f(t)^{\otimes k}\right) \le \left( \frac{k}{N} \right)^{\frac1q} \mathcal{W}^{[q]}_p\left(F^N(t)^s, f(t)^{\otimes N}\right).
\end{equation}

We now construct a (symmetric) coupling $\Pi_t\in \mathcal{C}(F^N(t)^s,f(t)^{\otimes N})\subset \mathcal{P}((\Omega_x \times \R^d_\xi)^{2N})$ in the following way: At the initial time we consider
\[
\Pi_0(Z^1_N,Z^2_N):=\delta_{(Z_N^2)}(Z_N^1)\otimes f(0)^{\otimes N}(Z_N^2),
\]
where $Z^i_N=(x^i_1,\xi^i_1,\dots,x_N^i,\xi_N^i)\in(\Omega_x\times\R_\xi^d)^N$, $X^i_N=(x_1^i, \dots, x^i_N) \in \Omega_x^{N}$, $\Xi_N^i=(\xi_1^i,\dots,\xi_N^i)\in\R^{dN}_\xi$ for $i=1,2$.
We then introduce  
\begin{equation}
\label{newcoupling}
\Pi_t(Z^1_N,Z^2_N):=\frac{1}{N!} \sum_{\sigma \in \mathfrak{S}_N}(\sigma \otimes \sigma)_*\left[\left( \Phi^N(t) \otimes \varphi_{f_0}(t)^{\otimes N}\right)_* \Pi_0(Z^1_N,Z^2_N)\right],
\end{equation}
where $\Phi^N(t)$ is the particle flow at time $t$ generated by the particle system in \eqref{sec0:maineq}, while $\varphi_{f_0}(t)$ is the Vlasov flow given in \eqref{sec1:mono}.
It is now easy to check that
$$
p^1_* \Pi_t=F^N(t)^s, \quad p_*^2\Pi_t=f(t)^{\otimes N},
$$
where $p^1$ and $p^2$ are the projections on the first and second copy of $(\Omega_x\times \R^d_\xi)^N\times (\Omega_x\times \R^d_\xi)^N$.

Hence, by \eqref{proof1:margest} and the definition of the coupling, we get
\begin{equation}
\label{proof1:estcoupl}
\begin{aligned}
\left(\frac{N}{k}\right)^{\frac{1}{q}}\mathcal{W}^{[q]}_p\left(F^{N}(t)^s_{N:k}, f(t)^{\otimes k}\right) 
&\le \left( \int \left(\text{d}^{(q)}_\Omega(X_N^1-X_N^2)+\| \Xi_N^1 -\Xi_N^2\|_{\ell_q}\right)^p d\Pi_t(Z^1_N,Z^2_N) \right)^{\frac{1}{p}}\\
& \le \left( \frac{1}{N!} \sum_{\sigma \in \mathfrak{S}_N} \int \| \sigma(\Xi_N^1(t)) -\sigma(\Xi_N^2(t))\|_{\ell_q}^p df(0)^{\otimes N}(Z_N^2)\right)^{\frac{1}{p}}
\end{aligned}
\end{equation}
where $\Xi_N^1(t):=\pi_\xi^{\otimes N}\left(\Phi^{N}(t)(Z_N^2)\right)$ and $\Xi_N^2(t):=\pi_\xi^{\otimes N} \left(\varphi_{f_0}(t)^{\otimes N}(Z_N^2)\right)$, with $\pi_\xi: \Omega_x \times \R^d_\xi \to \R^d_\xi$ the projection on the $\xi$ variable.

Let $\mu_N(t)$ be the empirical measure supported on the particle system at time $t$ with initial data $Z^2_N$. Similarly to the proof of part \eqref{sec0:thm1pt2} of Theorem \ref{sec0:thm1}, we have
\begin{align*}
\| \Xi_N^1(t) -\Xi_N^2(t)\|_{\ell_q} &\le \text{Lip}_{x,\xi}(G) \int_0^t\|\Xi_N^1(\tau) - \Xi_N^2(\tau) \|_{\ell_q }d \tau
+ N^{\frac{1}{q}} \text{Lip}_{x,\xi}(G) 
\int_0^t \mathcal{W}_1^{[q]}(\mu_N(\tau)^{\otimes m}, f(\tau)^{\otimes m}) d \tau\\
& \le \text{Lip}_{x,\xi}(G) \int_0^t \|\Xi_N^1(\tau) - \Xi_N^2(\tau) \|_{\ell_q }d \tau
+ (Nm)^{\frac{1}{q}} \text{Lip}_{x,\xi}(G) 
\int_0^t \mathcal{W}_1(\mu_N(\tau), f(\tau)) d \tau,
\end{align*}  
where we used that $\mathcal{W}^{[q]}_1\left(\mu_N(\tau)^{\otimes m}, f(\tau)^{\otimes m}\right) \le m^{\frac1q}\mathcal{W}_1(\mu_N(\tau), f(\tau))$ (see \cite[Lemma 14 (Appendix A.1.4)]{PT24}).
We now estimate
$$\mathcal{W}_1\left(\mu_N(\tau), f(\tau) \right) \le
 \mathcal{W}_1\left(\mu_N(\tau), \nu_N(\tau)\right)
 + \mathcal{W}_1\left( \nu_N(\tau), f(\tau) \right),
$$
where $\nu_N(\tau):=\varphi_{f_0}(\tau)_* \mu_N(0)$, with $\mu_N(0)$ the empirical measure supported on $Z^2_N$. Since the labels are stationary and they coincide at the initial time, by the Hölder's inequality,
\begin{align*}
 \mathcal{W}_1\left(\mu_N(\tau), \nu_N(\tau)\right)&\le \frac{1}{N} \sum_{i=1}^N |\xi^1_i(\tau)-\xi^2_i(\tau)|
 \le \frac{1}{N^{\frac{1}{q}}}\|\Xi_N^1(\tau)-\Xi_N^2(\tau)\|_{\ell_q},
\end{align*}
so that, by the Grönwall's lemma,
$$
\| \Xi_N^1(t) -\Xi_N^2(t)\|_{\ell_q}\le (Nm)^{\frac{1}{q}} \text{Lip}_{x,\xi}(G) e^{\text{Lip}_{x,\xi}(G)(1+{m}^{1/q}) t} \int_0^t \mathcal{W}_1\left(\nu_N(\tau), f(\tau) \right) d\tau.
$$ 
Since $\mathcal{W}_1 \le \mathcal{W}_p$, by \eqref{proof1:estcoupl} and the Minkowski's inequality we get
\begin{align*}
\mathcal{W}^{[q]}_p&(F^N(t)^s_{N:k}, f(t)^{\otimes k}) \\
&\le 
 (km)^{\frac{1}{q}} \text{Lip}_{x,\xi}(G) e^{\text{Lip}_{x,\xi}(G)(1+m^{\frac{1}{q}}) t}
\left( \int  \left( \int_0^t \mathcal{W}_p\left(\nu_N(\tau), f(\tau) \right) d\tau \right)^p df_0\right)^{\frac{1}{p}}\\
&\le
 (km)^{\frac{1}{q}} \text{Lip}_{x,\xi}(G) e^{\text{Lip}_{x,\xi}(G)(1+m^{1/q}) t}
\left( \int_0^t  \left( \int \mathcal{W}_p\left(\nu_N(\tau), f(\tau) \right)^p df_0 \right)^{\frac{1}{p}} d \tau \right).
\end{align*}
Since $\nu_N(\tau):=\varphi_{f_0}(\tau)* \mu_N(0)$ and $\nu^N(0)=\mu^N(0)$, we get that 
\[
 \int \mathcal{W}_p\left(\nu_N(\tau), f(\tau) \right)^p df(0) = \int \mathcal{W}_p\left(\mu_N(0), f(\tau) \right)^p df(\tau).
\]
 We now use a quantitative version of the law of large numbers given by the Fournier and Guillin's bound (see \cite{FG}) so that, for $z>0$ and $p\ge1$, we get 
 \begin{align*}
\int \mathcal{W}_p\left(\mu_N(0), f(\tau) \right)^p df(\tau)\le  
C \mathcal{M}_z(f(\tau))^{p/z}
\begin{cases}
N^{-\frac{1}{2}}+ N^{-(z-p)/z} \quad \text{if} \quad p>\frac{d}{2} \, \text{and} \, z\neq 2p\\
N^{-\frac{1}{2}}\log(1+N)+ N^{-(z-p)/z} \quad \text{if} \quad p=\frac{d}{2} \, \text{and} \, z\neq 2p \\
N^{-\frac{p}{d}}+ N^{-(z-p)/z} \quad \text{if} \quad p \in (0, \frac{d}{2})  \, \text{and} \, z\neq d/(d-p),
\end{cases}
\end{align*}
where $\mathcal{M}_z(f)=\int |x|^z df(x)$.\qed

\section{Proof of Theorem \ref{sec0:thm2}: Validation of the macroscopic equation and multiple-wise limit}
\label{sec4:sec}

\textbf{Proof of part \eqref{sec0:thm2pt1} and \eqref{sec0:thm2pt2} of Theorem \ref{sec0:thm2}:}
The proof of part \eqref{sec0:thm2pt1} follows by testing $f_{\text{mon}}(t,x,\xi)=\nu(x)\otimes\delta_{y(t,x)}(\xi)$ against a test function and noticing that
$$
\mathcal{X}^{[m]}[f_{\text{mon}}(t)](t,x,\xi)= \int_{\Omega^m} G^{(m)}(t,x,\xi,x_1,y(t,x_1),\dots, x_m, y(t,x_m)) d\nu^{\otimes m}(x_1,\dots,x_m).
$$

The existence and uniqueness of global solutions to \eqref{sec2:opinion} follows again by a contraction argument. Indeed, given $y_0 \in L^\infty(\Omega, \R^d)$, we define the operator
\[
\mathcal{F}: L^\infty\left([0,+\infty);L^\infty(\Omega_x,;\R^d_\xi) \right) \to L^\infty\left([0,+\infty);L^\infty(\Omega_x,;\R^d_\xi) \right)
\]
such that 
$$\mathcal{F}(\psi)(t,x):=y_0(x)+\int_0^t \int_{\Omega^m}G^{(m)}(s,x,\psi(s,x),x_1,\psi(s,x_1),\dots,x_m,\psi(s,x_m)) d\nu^{\otimes m}(x_1,\dots,x_m) ds.$$

By the Lipschitz regularity of $G^{(m)}$, it follows that, for any $C_m>0$,
{\small\begin{equation}
    \label{sec2:contraction}
\sup_{t\in[0,+\infty)} e^{-2C_m t}\|\mathcal{F}(\psi_1)(t)-\mathcal{F}(\psi_2)(t)\|_{L^\infty(\Omega;\R^d)}\le \frac{\text{Lip}_{\xi}(G^{(m)})(m+1)}{2C_m}\sup_{t\in[0,+\infty)}e^{-2C_m t}\|\psi_1-\psi_2\|_{L^\infty(\Omega,\R^d)}.
\end{equation}}
Choosing $C_m:=\text{Lip}_\xi(G^{(m)})(m+1)$, we obtain that $\mathcal{F}$ is a contraction in $L^\infty\left([0,+\infty);L^\infty(\Omega_x,;\R^d_\xi) \right)$ with respect to the weighted metric defined in \eqref{sec2:contraction}. 

Hence, there exists a unique $y^{(m)} \in L^\infty\left([0,+\infty);L^\infty(\Omega_x,;\R^d_\xi) \right)$ such that
{\small\begin{equation}
\label{sec2:intrep}
y^{(m)}(t,x)=y_0+\int_0^t \int_{\Omega^m}G^{(m)}(s,x,y^{(m)}(s,x),x_1,y^{(m)}(s,x_1),\dots,x_m,y^{(m)}(s,x_m)) d\nu^{\otimes m}(x_1,\dots,x_m) ds.
\end{equation}}
By the integral representation in \eqref{sec2:intrep} and the Lipschitz regularity of $G^{(m)}$, it follows that $y^{(m)}(\cdot,x)$ is also continuously differentiable in time and that \eqref{sec2:opinion} holds.\qed\\

\textbf{Proof of part \eqref{sec0:thm2pt3} of Theorem \ref{sec0:thm2}:}
In order to prove the last two parts of the Theorem, we first recall some results about the limit of symmetric functions of many variables in the following Proposition.
We refer to \cite[Theorem 2.1, Theorem 5.14]{CD} and \cite{HS55} for the proofs.
\begin{proposition}
\label{sec2:prop}
 Given $(\Omega, \text{d}_\Omega)$ a compact metric space and $m \in \N^*$, let $F^{(m)}:\Omega^m \to \R^d$ a symmetric function of the $m$ variables, i.e.
    \[
    F^{(m)}(x_1,\dots,x_m)=F^{(m)}(x_{\sigma(1)},\dots,x_{\sigma(m)})
    \]
    for any $\sigma$ permutation of $\{1,\dots,m\}$. 
The following results hold:
\begin{enumerate}[(i)]
    \item
   Assume that
    \begin{itemize}
    \item there exists $C_1>0$ independent of $m$ such that $\|F^{(m)}\|_{L^\infty(\Omega^m)}\le C_1$;
    \item there exists $C_2>0$ independent of $m$ such that 
    \[
     |F^{(m)}(X_m)-F^{(m)}(Y_m)| \le C_2 \mathcal{W}_1(\mu_{X_m},\mu_{Y_m}),
    \]
    for any $X_m=(x_1,\dots,x_m), Y_m=(y_1,\dots,y_m)\in\Omega^m$ and where $\mu_{X}=\frac1m\sum_{j=1}^m \delta_{x_j}$;
    \end{itemize}
    thene, there exists $F^\infty:\mathcal{P}(\Omega) \to \R^d$ and a subsequence $m_k \to +\infty$ as $k \to +\infty$ such that
    \begin{equation}
        \label{sec2:limsub}
    \lim_{k \to +\infty} \sup_{X_{m_k} \in \Omega^{m_k}}\left| F^{(m_k)}(X_{m_k}) - F^\infty(\mu_{X_{m_k}}) \right|=0.
    \end{equation}

    \item (Hewitt-Savage Theorem) Given $\nu \in \mathcal{P}(\Omega)$, there exist a subsequence $m_k \to +\infty$ as $k \to +\infty$ and $F^\infty:\mathcal{P}(\Omega)\to \R^d$ given by \eqref{sec2:limsub} such that
    \[
    \lim_{k \to +\infty} \int_{\Omega^{m_k}}F^{(m_k)}(x_1,\dots,x_{m_k}) d\nu^{\otimes m_k}(x_1,\dots,x_{m_k})=F^\infty(\nu).
    \]
    \end{enumerate}
\end{proposition}

We now study the limit of the macroscopic equation \eqref{sec2:opinion} as the order $m$ of the multiple-wise interaction goes to infinity. Let us write again the macroscopic equation 
\begin{equation}
\label{sec2:eq1}
\partial_t y^{(m)}(t,x) = \int_{\Omega_x^m} G^{(m)}(t,x,y^{(m)}(t,x),x_1,y^{(m)}(t,x_1),\dots,x_m,y^{(m)}(t,x_m)) d\nu^{\otimes m}(x_1,\dots,x_m),
\end{equation}
for $\nu \in \mathcal{P}(\Omega_x)$ and $m \in \N^*$. Notice that, since the measure $\nu^{\otimes m}$ is tensorized,
\begin{align*}
&\int_{\Omega_x^m} G^{(m)}(t,x,y^{(m)}(t,x),x_1,y^{(m)}(t,x_1),\dots,x_m,y^{(m)}(t,x_m)) d\nu^{\otimes m}(x_1,\dots,x_m)\\
&=\frac{1}{m!}\sum_{\sigma \in \mathfrak{S}_m}\int_{\Omega_x^m} G^{(m)}(t,x,y^{(m)}(t,x),x_{\sigma(1)},y^{(m)}(t,x_{\sigma(1)}),\dots,x_{\sigma(m)},y^{(m)}(t,x_{\sigma(m)})) d\nu^{\otimes m}(x_1,\dots,x_m),
\end{align*}
in particular, we can assume without loss of generality, that $G^{(m)}$ is symmetric with respect to the variables $\{x_i\}_{i=1}^m$.

Given $T>0$, by part \eqref{sec0:thm2pt2} of Theorem \ref{sec0:thm2}, we know there exists a unique solution $y^{(m)} \in C^1\left([0,T],L^\infty(\Omega_x;\R^d) \right)$ to \eqref{sec2:eq1} with initial datum $y_0 \in L^\infty(\Omega_x;\R^d)$. If moreover $y_0$ is Lipschitz in $\Omega_x$, thanks to the assumptions \eqref{sec0:ass1}-\eqref{sec0:ass2}, the sequence of solutions $\{y^{(m)}\}_{m=1}^\infty$ is uniformly bounded and Lipschitz equi-continuous: Indeed by \eqref{sec2:eq1}
\[
\sup_{t \in [0,T]}\|y^{(m)}(t,\cdot)\|_{L^\infty(\Omega_x)}\le \|y_0\|_{L^\infty(\Omega_x)}+T\|G^{(m)}\|_{L^\infty}\le \|y_0\|_{L^\infty(\Omega_x)}+TM,
\]
where $M$ is defined in \eqref{sec0:assM};
moreover
\begin{equation}
\label{sec3:est1}
\left|y^{(m)}(t_1,x_1)-y^{(m)}(t_2,x_2)\right|\le M|t_1-t_2| + e^{\text{Lip}_{x,\xi}G T}(\text{Lip}_xy_0+T\text{Lip}_{x,\xi}G)|x_1-x_2|.
\end{equation}
By the Arzéla-Ascoli Theorem, we deduce there exists $y^\infty \in C([0,T],C(\Omega_x;\R^d))$ Lipschitz-continuous with the same constants obtained in \eqref{sec3:est1} such that, up to a subsequence, $y^{(m)}$ converges to $y^\infty$ in $L^\infty$.

We now want to show that $y^\infty$ satisfies \eqref{sec1:limitm}. To do this, we define 
\[
\widetilde G^{(m)}_{t,x}[y^{\infty}](x_1,\dots,x_m):=G^{(m)}(t,x,y^\infty(t,x), x_1, y^{\infty}(t,x_1),\dots,x_m,y^\infty(t,x_m))
\]
and we notice $\tilde G_{t,x}^{(m)}[y^\infty]$ is symmetric in the variables $x_i\in\Omega$. Thanks to \eqref{sec0:assM} and the Lipschitz-continuity of $y^\infty$ given by \eqref{sec3:est1}, $\widetilde G^{(m)}_{t,x}[y^\infty]$ satisfies the assumptions of part (i) of Proposition \ref{sec2:prop}. 

Hence
there exists a continuous $G^\infty_{t,x}[y^\infty]:\mathcal{P}(\Omega) \to \R^d $ such that \eqref{sec2:limsub} holds.
By the triangle inequality on $G^{(m)}$ and by the Hewitt-Savage Theorem in Proposition \ref{sec2:prop}, 
we get
\begin{align*}
\lim_{m \to +\infty}\int_{\Omega^m} G^{(m)}(t,x,y(t,x),&x_1,y(t,x_1),\dots,x_m,y(t,x_m)) d\nu^{\otimes m}(x_1,\dots,x_m)=G^\infty_{t,x}[y^\infty](\nu).\\
\end{align*}

\textbf{Proof of part \eqref{sec0:thm2pt4} of Theorem \ref{sec0:thm2}:} Given $\phi \in \mathcal{C}(\Omega_x)$, by triangle inequality we have
\begin{align}
    \left|\int_{\Omega_x} \phi(x)\left(Y_{N:1}^{m}(t,dx)-y^\infty(t,x)\nu(dx) \right)\right| &\le \left|\int_{\Omega_x} \phi(x)\left(Y_{N:1}^{m}(t,dx)-y^{(m)}(t,x)\nu(dx) \right)\right|\\
    &+ \left|\int_{\Omega_x} \phi(x)\left(y^{(m)}(t)\nu(dx)-y^\infty(t,x)\nu(dx) \right)\right|\equiv T_1+T_2.
\end{align}
By the characterization of weak convergence through the Wasserstein distance $ \mathcal{W}^{[1]}_1$, we get
\begin{align}
    T_1=\Bigl|\int_{\Omega_x\times \R^d} \phi(x)\xi&\Bigl(dF^{N,(m)}(t)^s_{N:1}(x,\xi)-\delta(\xi-y^{(m)}(t,x))\otimes d\nu(x) \Bigr)\Bigr|\\
    &\le \|\phi\|_{L^\infty}
    \sup_{t \in [0,T]}\mathcal{W}^{[1]}_1(F^{N,(m)}(t)^s_{N:1}, \delta(\xi-y^{(m)}(t,x))\nu(dx)),
\end{align}
while, since $\Omega_x$ is compact,
\begin{equation}
    T_2 \le \|\phi\|_{L^\infty}\sup_{t \in [0,T]}\|y^{(m)}(t,\cdot)-y^\infty(t,\cdot)\|_{L^\infty(\Omega)}.
\end{equation}
By taking $m\le m(N)$ as specified in \eqref{sec2:m(N)}, the conclusion follows by the results in \eqref{sec0:result2} and \eqref{sec3:multilim}.
\qed

\bibliographystyle{plain}
\bibliography{bibliography}

\begin{thebibliography}{1}

\bibitem{BATTISTON20201}
F.~Battiston, G.~Cencetti, I.~Iacopini, V.~Latora, M.~Lucas, A.~Patania, J.-G.
  Young, and G.~Petri.
\newblock Networks beyond pairwise interactions: Structure and dynamics.
\newblock {\em Physics Reports}, 874:1--92, 2020.

\bibitem{CD}
P.~Cardaliaguet.
\newblock Notes on mean field games.
\newblock {\em https://www.ceremade.dauphine.fr/~cardaliaguet/MFG20130420.pdf},
  2013.

\bibitem{Dobr}
R.~L. Dobrušin.
\newblock Vlasov equations.
\newblock {\em Funktsional. Anal. i Prilozhen}, 13:48--58, 1979.

\bibitem{FG}
N.~Fournier and A.~Guillin.
\newblock On the rate of convergence in {W}asserstein distance of the empirical
  measure.
\newblock {\em Probab. Theory Related Fields}, 162(3-4):707--738, 2015.

\bibitem{GMP16}
F.~Golse, C.~Mouhot, and T.~Paul.
\newblock On the mean field and classical limits of quantum mechanics.
\newblock {\em Comm. Math. Phys.}, 343(1):165--205, 2016.

\bibitem{HS55}
E.~Hewitt and L.~J. Savage.
\newblock Symmetric measures on {C}artesian products.
\newblock {\em Trans. Amer. Math. Soc.}, 80:470--501, 1955.

\bibitem{RUSSO24}
F.~Malizia, A.~Corso, and L.V. Gambuzza~et al.
\newblock Reconstructing higher-order interactions in coupled dynamical
  systems.
\newblock {\em Nat. Commun.}, 5184(15), 2024.

\bibitem{PT24}
T.~Paul and E.~Trélat.
\newblock From microscopic to macroscopic scale equations: mean field,
  hydrodynamic and graph limits.
\newblock {\em arXiv:2209.08832}, 2024.

\end{thebibliography}

\end{document}